\documentclass[12pt]{article}

\usepackage[latin1]{inputenc}
\usepackage{amsmath,amssymb}
\usepackage{latexsym}

\newtheorem{theorem}{Theorem}[section]
\newtheorem{corollary}[theorem]{Corollary}
\newtheorem{lemma}[theorem]{Lemma}
\newtheorem{proposition}[theorem]{Proposition}
\newtheorem{definition}[theorem]{Definition}

\numberwithin{equation}{section}

\def\proof{{\medskip\noindent {\bf Proof. }}}
\def\longproof#1{{\medskip\noindent {\bf Proof #1.}}}
\def\qed{{\hfill $\square$ \bigskip}}

\def\square{{\vcenter{\vbox{\hrule height.3pt
        \hbox{\vrule width.3pt height5pt \kern5pt
           \vrule width.3pt}
        \hrule height.3pt}}}}

 \def\sB {{\cal B}} 
  \def\sF {{\cal F}}
  \def\sI {{\cal I}}
 \def\sK {{\cal K}} \def\sL {{\cal L}}
\def\sM {{\cal M}}  \def\sO {{\cal O}}
\def\sP {{\cal P}}  
\def\sS {{\cal S}} \def\sT {{\cal T}}

\def\wt{\widetilde}

\def\ol{\overline}

\def\P{{\mathbb P}}

\def\bee{\begin{equation}}
\def\bet{\begin{theorem}}
\def\bep{\begin{proposition}}
\def\bel{\begin{lemma}}
\def\bec{\begin{corollary}}
\def\bed{\begin{definition}}
\def\eee{\end{equation}}
\def\eet{\end{theorem}}
\def\eep{\end{proposition}}
\def\eel{\end{lemma}}
\def\eec{\end{corollary}}
\def\eed{\end{definition}}

\def\R{{\mathbb R}}

\def\P{{\mathbb P}}

\def\F{{\cal F}}

\def\proof{{\medskip\noindent {\bf Proof. }}}
\def\longproof#1{{\medskip\noindent {\bf Proof #1.}}}
\def\qed{{\hfill $\square$ \bigskip}}

\def\eps{\varepsilon}

 \def\qq {\qquad}

\def\wt{\widetilde}
\def\ol{\overline}

\def\ni{\noindent }

\def\square{{\vcenter{\vbox{\hrule height.3pt
        \hbox{\vrule width.3pt height5pt \kern5pt
           \vrule width.3pt}
        \hrule height.3pt}}}}

\def\tlint{{- \kern-0.85em \int \kern-0.2em}}  % for textstyle
\def\dlint{{- \kern-1.05em \int \kern-0.4em}}  % for displays

 \def\sB {{\cal B}} 
  \def\sF {{\cal F}}
  \def\sI {{\cal I}}
 \def\sK {{\cal K}} \def\sL {{\cal L}}
\def\sM {{\cal M}}  \def\sO {{\cal O}}
\def\sP {{\cal P}}  
\def\sS {{\cal S}} \def\sT {{\cal T}}

\begin{document}

\title{The measurability of hitting times\\ (corrected version)}

\author{Richard F. Bass\footnote{Research partially supported by NSF grant
DMS-0901505.}}

\maketitle

\begin{abstract}  
\noindent Under very general conditions the hitting time of a set by
a stochastic process is a stopping time. We give a new simple proof of
this fact. The section theorems for
optional and predictable sets are easy corollaries of the proof.
\end{abstract}

\section{Introduction}

A fundamental theorem in the foundations of stochastic processes is
the one that says that, under very general  conditions,  the first time a stochastic process  enters
a set is a stopping time.
The proof uses capacities, analytic sets, and Choquet's capacibility theorem,
and is considered hard. To the best of our knowledge, no more than
a handful of books have an exposition that starts with the definition of capacity
and proceeds to the hitting time theorem. (One that does is \cite{DM}.)

The purpose of this paper is to give a short and elementary proof of this
theorem. The proof is simple enough that it could
easily be included in a first year graduate course in probability. 

In Section \ref{S:D} we give a proof of the debut theorem, from which 
the measurability theorem follows. 
As easy corollaries we obtain the section theorems for optional
and predictable sets. This argument is given in Section \ref{S:S}.

Note that this paper is a version of \cite{bass-choq1}, revised to
take into account the corrections in \cite{bass-choq2}.

\section{The debut theorem}\label{S:D}

Suppose $(\Omega, \F, \P)$ is a probability space.
The outer probability $\P^*$ associated with $\P$ is given by
$$\P^*(A)=\inf\{\P(B): A\subset B, B\in \sF\}.$$
A set $A$ is a $\P$-null set if $\P^*(A)=0$.
Suppose $\{\sF_t\}$ is
a filtration satisfying the usual conditions:  $\cap_{\eps>0} \F_{t+\eps}=\F_t$
for all $t\geq 0$, and each $\F_t$ contains every $\P$-null set.
 Let $\pi:[0,\infty)\times \Omega\to \Omega$
be defined by
$\pi(t,\omega)=\omega$.

Recall that a random variable taking values in $[0,\infty]$ is a stopping time if
$(T\leq t)\in \sF_t$ for all $t$; we allow our stopping times to take
the value infinity. Since the filtration satisfies the usual
conditions, $T$ will be a stopping time  if $(T<t)\in \sF_t$ for all $t$.
If $T_i$ is a finite collection or countable collection of stopping times,
then $\sup_i T_i$ and $\inf_i T_i$ are also stopping times.

Given a topological space $\sS$, the Borel $\sigma$-field is the
one generated by the open sets. 
Let
$\sB[0,t]$ denote the Borel $\sigma$-field
on $[0,t]$ and  $\sB[0,t]\times \sF_t$ the product $\sigma$-field.
A process $X$ taking values in a topological space 
$\sS$ is progressively measurable if for each $t$ the 
map $(s,\omega)\to X_s(\omega)$ from $[0,t]\times \Omega$ to $\sS$
is measurable with respect to
$\sB[0,t]\times\sF_t$, that is, the inverse image of Borel subsets of $\sS$
are elements of
 $\sB[0,t]\times \sF_t$.
If  the paths of $X$ are right continuous,
then $X$ is easily seen to be progressively measurable. The same is
true if $X$ has left continuous paths. A subset of $[0,\infty)\times \Omega$
is progressively measurable if its indicator is a progressively measurable
process.

If $E\subset [0,\infty)\times \Omega$, let $D_E=\inf\{t\geq 0: (t,\omega)\in E\}$, the debut of $E$.
%If $T$ is a stopping time, let $[T,T]=\{(t,\omega): t=T(\omega)<\infty\}$.
We will prove 

\bet\label{GT1} If $E$ is a progressively measurable set, then $D_E$ is
a stopping time.
\eet

Fix $t$. Let $\sK^0(t)$ be  the collection of subsets of $[0,t]\times \Omega$
of the form $K\times C$, where $K$ is a compact subset of $[0,t]$ and
$C\in \sF_t$. Let $\sK(t)$ be the collection of finite unions of sets in 
$\sK^0(t)$ and let $\sK_\delta(t)$ be the collection of countable intersections
 of
sets in $\sK(t)$.
We say $A\in \sB[0,t]\times \sF_t$ is $t$-approximable if given $\eps>0$, there exists
$B\in \sK_\delta(t)$ with $B\subset A$ and 
\bee\label{GLE1}
\P^*(\pi(A))\leq \P^*(\pi(B))+\eps.
\eee

\bel\label{GL1} If $B\in \sK_\delta(t)$, then $\pi(B)\in \sF_t$.
If $B_n\in \sK_\delta(t)$ and $B_n\downarrow B$, then $\pi(B)=\cap_n \pi(B_n)$.
\eel

\ni The hypothesis that the $B_n$ be in $\sK_\delta(t)$ is important. For example,
if $B_n=[1-(1/n),1)\times \Omega$, then $\pi(B_n)=\Omega$ but $\pi(\cap_n B_n)
=\emptyset$. This is why the proof given in \cite[Lemma 6.18]{Elliott}
is incorrect.

\proof If $B=K\times C$, where $K$ is a nonempty subset of $[0,t]$ and 
$C\in \sF_t$, then $\pi(B)=C\in \sF_t$. Therefore $\pi(B)\in \sF_t$
if $B\in \sK^0(t)$. If $B=\cup_{i=1}^m A_i$ with $A_i\in \sK^0(t)$, then
$\pi(B)=\cup_{i=1}^m \pi(A_i)\in \sF_t$. 

For each $\omega$ and each set $C$, let 
\bee\label{GLE2}
S(C)(\omega)=\{s\leq t: (s,\omega)\in C\}.
\eee 
If $B\in \sK_\delta(t)$ and $B_n\downarrow B$ with $B_n\in \sK(t)$ for
each $t$, then $S(B_n)(\omega)\downarrow S(B)(\omega)$, so $S(B)(\omega)$
is compact. 

Now suppose $B\in \sK_\delta(t)$ and take $B_n\downarrow B$
with $B_n\in \sK_\delta(t)$. 
$S(B_n)(\omega)$ is a compact subset of $[0,t]$
for each $n$ and $S(B_n)(\omega)\downarrow S(B)(\omega)$. One possibility is that
$\cap_n S(B_n)(\omega)\ne \emptyset$; in this case, if  $s\in \cap_n S(B_n)(\omega)$,
then $(s,\omega)\in B_n$ for each $n$, and so $(s,\omega)\in B$. Therefore
$\omega\in \pi(B_n)$ for each $n$ and $\omega\in \pi(B)$.
The other possibility is that $\cap_n S(B_n)(\omega)=\emptyset$. Since
the sequence $S(B_n)(\omega)$ is a decreasing sequence of compact sets, 
  $S(B_n)(\omega)=\emptyset$ for some $n$,
for otherwise $\{S(B_n)(\omega)^c\}$ would be an open cover of $[0,t]$
with no finite subcover.
Therefore $\omega\notin \pi(B_n)$ and $\omega\notin \pi(B)$. We conclude
that $\pi(B)=\cap_n \pi(B_n)$. 

Finally, suppose $B\in \sK_\delta(t)$ and $B_n\downarrow B$ with $B_n\in \sK(t)$.
Then $\pi(B)=\cap_n \pi(B_n)\in \sF_t$.  
\qed

\bep\label{GP2} Suppose $A$ is $t$-approximable. Then $\pi(A)\in \sF_t$.
Moreover, given $\eps>0$ there exists $B\in \sK_\delta(t)$ such that $\P
(\pi(A)\setminus \pi(B))<\eps$.
\eep

\proof Choose $A_n\in \sK_\delta(t)$ with  $A_n\subset A$ and
$\P(\pi(A_n))\to \P^*(\pi(A))$. Let $B_n=A_1\cup \cdots \cup A_n$ and
let $B=\cup_n B_n$. Then
$B_n\in \sK_\delta(t)$, $B_n\uparrow B$, and $\P(\pi(B_n))\geq \P(\pi(A_n))\to 
\P^*(\pi(A))$. It follows that
$\pi(B_n)\uparrow \pi(B)$, and so $\pi(B)\in \sF_t$ and 
$$\P(\pi(B))=\lim \P(\pi(B_n))=\P^*(\pi(A)).$$ 

For each $n$, there exists $C_n\in \sF$ such that
$\pi(A)\subset C_n$ and $\P(C_n)\leq \P^*(\pi(A))+1/n$. Setting
$C=\cap_n C_n$, we have $\pi(A)\subset C$ and $\P^*(\pi(A))=\P(C)$.
Therefore
$\pi(B)\subset \pi(A)\subset C$ and $\P(\pi(B))=\P^*(\pi(A))=\P(C)$.
This implies that $\pi(A)\setminus \pi(B)$ is a $\P$-null set,  and by
the completeness assumption, $\pi(A)=(\pi(A)\setminus \pi(B))\cup \pi(B)\in \sF_t$. 
Finally,
 $$\lim_n\P(\pi(A)\setminus \pi(B_n))=\P(\pi(A))\setminus \pi(B))=0.$$
\qed

The following lemma is well known; see, e.g., \cite[p.~94]{Bichteler}.

\bel\label{A.1} (a) If $A\subset \Omega$, there exists $C\in \sF$
such that $A\subset C$ and $\P^*(A)=\P(C)$.

(b) Suppose $A_n\uparrow A$. Then $\P^*(A)=\lim_{n\to \infty} \P^*(A_n)$.
\eel

\proof (a) By the definition of $\P^*(A)$, for each $n$ there exists
$C_n\in \sF$ such that $A\subset C_n$ and $\P(C_n)\leq \P^*(A)+(1/n)$.
Setting $C=\cap_n C_n$, we have $A\subset C$, $C\in \sF$, and
$\P(C)\leq \P(C_n)\leq \P^*(A)+(1/n)$ for each $n$, hence $\P(C)\leq \P^*(A)$.

(b) 
Choose $C_n\in \sF$ with $A_n\subset C_n$ and $\P^*(A_n)=\P(C_n)$.
Let $D_n=\cap_{k\ge n} C_k$ and $D=\cup_n D_n$. We see that
$D_n\uparrow D$, $D\in \sF$, and $A\subset D$. Then
\bee\tag*{$\square$}\P^*(A)\ge \sup_n \P^*(A_n)=\sup_n \P(C_n)\geq \sup_n \P(D_n)
=\P(D)\geq \P^*(A).\eee

Let $\sT_t=[0,t]\times \Omega$. Given a compact Hausdorff space $X$,
let $\rho^X: X\times \sT_t\to \sT_t$ be defined by $\rho^X(x,(s,\omega))
=(s,\omega)$. Let 
$$\sL_0(X)=\{A\times B: A\subset X, A\mbox{ compact},
B\in \sK(t)\},$$
$\sL_1(X)$ the class of finite unions of  sets in $\sL_0(X)$, and 
$\sL(X)$ the class of intersections of countable decreasing sequences
in $\sL_1(X)$. Let $\sL_\sigma(X)$ be the class of unions of countable increasing
sequences of sets in $\sL(X)$ and $\sL_{\sigma\delta}(X)$ the class
of intersections of countable decreasing sequences of sets in $\sL_\sigma(X)$.

\bel\label{A.2} If $A\in \sB[0,t]\times \sF_t$, there exists a compact
Hausdorff space $X$ and $B\in \sL_{\sigma\delta}(X)$ such that
$A=\rho^X(B)$.
\eel

\proof If $A\in \sK(t)$, we take $X=[0,1]$, the unit interval with the
usual topology and $B=X\times A$. Thus the collection $\sM$ of 
subsets of $\sB[0,t]\times \sF_t$
for which the lemma is satisfied contains $\sK(t)$. We will
show that $\sM$ is a monotone class.

Suppose $A_n\in \sM$ with $A_n\downarrow A$. There exist compact Hausdorff
spaces $X_n$ and sets  $B_n\in \sL_{\sigma\delta}(X_n)$ such that $A_n=\rho^{X_n}(B_n)$.
Let $X=\prod_{n=1}^\infty X_n$ be furnished with the product topology. Let
$\tau_n: X\times \sT_t\to X_n\times \sT_t$ be defined by $\tau_n(x,(s,\omega))
=(x_n,(s,\omega))$ if $x=(x_1,x_2, \ldots)$. Let $C_n=\tau_n^{-1}(B_n)$
and let $C=\cap_n C_n$. It is easy to check that $\sL(X)$ is closed under
the operations of finite unions and intersections, from which it follows
that $C\in \sL_{\sigma\delta}(X)$. If $(s,\omega)\in A$, then for each $n$ there exists $x_n\in X_n$ such that $(x_n,(s,\omega))\in B_n$. Note that
$((x_1,x_2, \ldots),(s,\omega))\in C$ and therefore $(s,\omega)\in \rho^X(C)$.
It is straightforward that $\rho^X(C)\subset A$, and we conclude 
$A\in \sM$.

Now suppose $A_n\in \sM$ with $A_n\uparrow A$. Let $X_n$ and $B_n$ be as before.
Let $X'=\cup_{n=1}^\infty (X_n\times \{n\})$ with the topology generated by
$\{G\times \{n\}: G\mbox{ open in } X_n\}$. Let $X$ be the one point
compactification of $X'$. We can write $B_n=\cap_m B_{nm}$ with
$B_{nm}\in \sL_\sigma(X_n)$. Let
$$C_{nm}=\{((x,n),(s,\omega))\in X\times \sT_t: x\in X_n, (x,(s,\omega))\in B_{nm}\},$$
$C_n=\cap_m C_{nm}$, and $C=\cup_n C_n$.
Then $C_{nm}\in \sL_\sigma(X)$ and  so $C_n\in \sL_{\sigma \delta}(X)$.

If $((x,p),(s,\omega))\in \cap_m \cup_n C_{nm}$, then 
for each $m$ there exists $n_m$ such that $((x,p),(s,\omega))\in C_{n_mm}$.
 This is
only possible if $n_m=p$ for each $m$. Thus $((x,p), (s, \omega))\in \cap_m  C_{pm}=C_p\subset C$.
The other inclusion is easier and we thus obtain $C=\cap_m\cup_n C_{nm}$,
which implies $C\in \sL_{\sigma\delta}(X)$. We check that
$A=\rho^X(C)$ along the same lines, and therefore $A\in \sM$.

If $\sI^0(t)$ is the collection of sets of the form $[a,b)\times C$, where
$a<b\leq t$ and $C\in \sF_t$, and $\sI(t)$ is the collection of finite
unions of sets in $\sI^0(t)$, then $\sI(t)$ is an algebra of sets. We
note that $\sI(t)$ generates the $\sigma$-field $\sB[0,t]\times \sF_t$. A set
in $\sI^0(t)$ of the form $[a,b)\times C$ is the union of sets in $\sK^0(t)$
of the form $[a, b-(1/m)]\times C$,  and it
follows that every set in $\sI(t)$ is the increasing union of sets
in $\sK(t)$. Since $\sM$ is a monotone
class containing $\sK(t)$, then $\sM$ contains $\sI(t)$.
By the monotone class theorem, $\sM=\sB[0,t]\times \sF_t$.
\qed

The works of Suslin and Lusin present a different approach to the
idea of representing Borel sets as projections; see, e.g., 
\cite[p.~88]{Fonseca} or \cite[p.~284]{Cohn}.

\bel\label{LA3} If $A\in \sB[0,t]\times \sF_t$, then $A$ is
$t$-approximable.
\eel

\proof
We first prove that if $H\in \sL(X)$, then $\rho^X(H)\in \sK_\delta$. If $H\in \sL_1(X)$,
this is clear. Suppose that $H_n\downarrow H$ with each $H_n\in \sL_1(X)$.
If $(s,\omega)\in \cap_n \rho^X(H_n)$, there exist
$x_n\in X$ such that $(x_n,(s,\omega))\in H_n$. Then there exists a subsequence
such that $x_{n_k}\to x_\infty$ by the compactness of $X$. Now $(x_{n_k},(s,\omega))\in H_{n_k}
\subset H_m$ for $n_k$ larger than $m$.  For fixed $\omega$, $\{(x,s): (x,(s,\omega))\in H_m\}$
is compact, so $(x_\infty,(s,\omega))\in H_m$ for all $m$. This implies
$(x_\infty,(s,\omega))\in H$. The other inclusion is easier  and therefore $\cap_n \rho^X(H_n)=\rho^X(H)$.
Since $\rho^X(H_n)\in \sK_\delta(t)$, then $\rho^X(H)\in \sK_\delta(t)$.
We also observe that for fixed $\omega$, $\{(x,s):(x,(s,\omega))\in H\}$
is compact.

Now suppose $A\in \sB[0,t]\times \sF_t$. Then by Lemma
\ref{A.2} there exists a compact Hausdorff  space $X$ and $B\in \sL_{\sigma \delta}(X)$ such that $A=\rho^X(B)$. We can write
$B=\cap_n B_n$ and $B_n=\cup_m B_{nm}$ with $B_n\downarrow B$, $B_{nm}\uparrow B_n$, and $B_{nm}\in \sL(X)$.

Let $a=\P^*(\pi(A))=\P^*(\pi\circ \rho^X(B))$ and let $\eps>0$. 
By Lemma \ref{A.1},
$$\P^*(\pi\circ \rho^X(B\cap B_{1m}))\uparrow \P^*(\pi\circ \rho^X(B\cap B_1))
=\P^*(\pi\circ \rho^X(B))=a.$$
Take $m$ large enough so that $\P^*(\pi\circ \rho^X(B\cap B_{1m}))>a-\eps$,
let $C_1=B_{1m}$, and $D_1=B\cap C_1$.

We proceed by induction. Suppose we are given sets $C_1, \ldots, C_{n-1}$ and 
sets $ D_1, \ldots, D_{n-1}$ with $D_{n-1}=B\cap (\cap_{i=1}^{n-1} C_i)$, $\P^*(\pi
\circ \rho^X(D_{n-1}))>a-\eps$,  and each $C_i=B_{im_i}$ for
some $m_i$. Since $D_{n-1}\subset B\subset B_n$, by Lemma \ref{A.1}
$$\P^*(\pi\circ \rho^X(D_{n-1}\cap B_{nm}))
\uparrow \P^*(\pi\circ \rho^X(D_{n-1}\cap B_n))
=\P^*(\pi\circ \rho^X(D_{n-1})).$$
We can take $m$ large enough so that
$\P^*(\pi\circ \rho^X(D_{n-1}\cap B_{nm}))>a-\eps$, let $C_n=B_{nm}$, and $D_n=D_{n-1}\cap C_n$.

If we let $G_n=C_1\cap \cdots \cap C_n$ and $G=\cap_n G_n=\cap_n C_n$, then
 each $G_n$ is in $\sL(X)$, hence $G\in \sL(X)$. Since $C_n\subset B_n$, then
$G\subset \cap_n B_n=B$. 
Each $G_n\in \sL(X)$ and so by the first paragraph of this proof, for each
fixed $\omega$ and $n$, $\{(x,s): (x,(s,\omega))\in G_n\}$
is compact. Hence by a proof very similar to that of Lemma 2.2, 
$\pi\circ \rho^X(G_n)\downarrow \pi\circ\rho^X(G)$.
 Using the first paragraph of this proof and Lemma 2.2,
we see  that $$\P(\pi\circ \rho^X(G))
=\lim \P(\pi\circ \rho^X(G_n))\geq \lim \P^*(\pi\circ \rho^X(D_n))\geq a-\eps.$$

Using the first paragraph of this proof once again, we see that $A$ is $t$-approxi\-mable. 
\qed

%\newpage

\longproof{of Theorem 2.1} 
Let $E$ be a progressively measurable  set and let $A_u=E\cap ([0,u]\times \Omega)$.
By Lemma \ref{LA3}, $A_u$ is $u$-approximable.
By Proposition 2.3, $\pi(A_u)\in \sF_u$.
Now fix $t$. If $\omega\in (D_E\le t)$, we see that  $\omega\in \pi(A_u)$ for all $u>t$.
Conversely, if $\omega\in \pi(A_u)$ for all $u>t$, note $\omega\in (D_E\le t)$. 
If $u_1<u_2$, then $A_{u_1}\subset A_{u_2}$ and hence $\pi(A_{u_1})\subset 
\pi(A_{u_2})$. Therefore                
$$(D_E\le t)=\cap_{u>t} \pi(A_u)\in \cap_{u>t} \sF_u=\sF_t.$$
Because $t$ was arbitrary,  we conclude $D_E$ is a stopping time.
\qed

If $B$ is a Borel subset of a topological space $\sS$, let
$$U_B=\inf\{t\geq 0: X_t\in B\}$$
and $$T_B=\inf\{t>0: X_t\in B\},$$
the first entry time and first hitting time of $B$, resp.

Here is the measurability theorem.

\bet\label{GC2} If $X$ is a progressively measurable process  taking values in 
$\sS$  and $B$ is a Borel subset of
$\sS$, then $U_B$ and $T_B$
are stopping times.
\eet

\proof Since $B$ is a Borel subset of $\sS$ and 
 $X$ is progressively measurable, then $1_B(X_t)$ is
also progressively measurable. $U_B$ is then the debut of the
set $E=\{(s,\omega):1_B(X_s(\omega))=1\}$, and therefore is a stopping time.

If we let $Y_t^\delta=X_{t+\delta}$ and $U^\delta_B=\inf\{t\geq 0:
Y_t^\delta\in B\}$, then by the above, $U^{\delta}_B$ is a stopping
time with respect to the filtration $\{\sF^\delta_t\}$, where $\sF_t^\delta=
\sF_{t+\delta}$. It follows that  
$\delta+U^\delta_B$ is a stopping time with respect to the filtration
$\{\sF_t\}$. Since $(1/m)+U_B^{1/m}\downarrow T_B$, then
$T_B$ is a stopping time with respect to $\{\sF_t\}$ as well.
\qed

We remark that in the theory of Markov processes, the notion of completion
of a $\sigma$-field is a bit different. In that case, we suppose
that $\sF_t$ contains all sets $N$ such that $\P^\mu(N)=0$ for every starting
measure $\mu$. 
The proof in Proposition \ref{GP2} shows that
$$(\P^\mu)^*(\pi(A)\setminus \pi(B))=0$$
for every starting measure $\mu$, so $\pi(A)\setminus\pi(B)$ is a $\P^\mu$-null
set for every starting measure $\mu$. Therefore $\pi(A)=\pi(B)\cup
(\pi(A)\setminus \pi(B))\in \sF_t$. With this modification, the rest of the
proof of Theorem \ref{GT1} goes through in the Markov process context.

\section{The section theorems}\label{S:S}

 Let $(\Omega, \sF, \P)$ be
a probability space and let $\{\sF_t\}$ be a filtration satisfying the usual conditions. 
The optional $\sigma$-field $\sO$ is the $\sigma$-field of subsets
of $[0,\infty)\times \Omega$ generated by the set of maps
$X: [0,\infty)\times \Omega\to \R$ where $X$ is bounded, adapted to
the filtration $\{\sF_t\}$, and has right continuous paths. The
predictable $\sigma$-field $\sP$ is the $\sigma$-field of subsets
of $[0,\infty)\times \Omega$ generated by the set of maps
$X: [0,\infty)\times \Omega\to \R$ where $X$ is bounded, adapted to
the filtration $\{\sF_t\}$, and has left continuous paths.

Given a stopping time $T$, we define $[T,T]=\{(t,\omega): t=T(\omega)<\infty\}$.
A stopping time is predictable if there exist stopping times $T_1, T_2, \ldots$
with $T_1\leq T_2\leq \cdots$, $T_n\uparrow T$, and  on the event $(T>0)$,
$T_n<T$ for all $n$.
We say the stopping times $T_n$ predict $T$. 
If $T$ is a predictable stopping time and $S=T$ a.s., we also call
$S$ a predictable stopping time.

The optional section theorem is the following.

\bet\label{SST1}
If $E$ is an optional set and $\eps>0$, there exists a stopping time $T$ such that
$[T,T]\subset E$ and $\P(\pi(E)) \leq \P(T<\infty)+\eps$.
\eet

The statement of the predictable section theorem is very similar.

\bet\label{SST2}
If $E$ is a predictable set and $\eps>0$, there exists a predictable stopping time $T$ such that
$[T,T]\subset E$ and $\P(\pi(E)) \leq \P(T<\infty)+\eps$.
\eet

First we prove the following lemma.

\bel\label{SSL1}
(1) $\sO$ is generated by the collection of processes $1_C(\omega)1_{[a,b)}(t)$
where $C\in \sF_a$.

(2) $\sP$ is generated by the collection of processes $1_C(\omega)1_{[b,c)}(t)$
where $C\in \sF_a$ and $a<b<c$.
\eel

\proof (1) First of all, $1_C(\omega)1_{[a,b)}(t)$ is a bounded right
continuous adapted process, so it is optional. 

Let $\sO'$ be the $\sigma$-field on $[0,\infty)\times \Omega$ generated by
the collection of processes $1_C(\omega)1_{[a,b)}(t)$, where $C\in \sF_a$.
Letting $b\to \infty$, $\sO'$ includes sets of the form $[a,\infty)\times C$ with
$C\in\sF_a$.

Let $X_t$ be a right continuous, bounded, and adapted process and let $\eps>0$.
 Let $U_0=0$ and define
$U_{i+1}=\inf\{t>U_i: |X_t-X_{U_i}|>\eps\}$. Since $(U_1<t)=\cup(|X_q-X_0|>\eps)$,
where the union is over all rational $q$ less than $t$, $U_1$ is a stopping 
time, 
and an analogous argument shows that each $U_i$ is also a stopping time.
If $S$ and $T$ are stopping times, let
$1_{[S,T)}=\{(t,\omega)\in [0,\infty)\times \Omega: S(\omega\leq t<T(\omega)\}$.
If we set 
$$X_t^\eps(\omega)=\sum_{i=0}^\infty X_{U_i}(\omega) 1_{[U_i,U_{i+1})}(t),$$
then $\sup_{t\geq 0}|X_t-X_t^\eps|\leq \eps$. Therefore we can approximate
$X$ by processes of the form
$$\sum_{i=0}^\infty X_{U_i} 1_{[U_i,\infty)}-\sum_{i=0}^\infty X_{U_i} 1_{[U_{i+1},\infty)}.$$
It therefore  suffices to show that if $V$ is a stopping time and $A\in \sF_V$,
then $1_A(\omega)1_{[V,\infty)}(t)$ is $\sO'$
measurable.

Letting $V_n=(k+1)/2^n$ when $k/2^n\leq V<(k+1)/2^n$,
\begin{align*}
1_A(\omega)1_{[V(\omega),\infty)}(t)&=\lim_{n\to \infty}
1_A(\omega)1_{[V_n(\omega),\infty)}(t)\\
&=\lim_{n\to \infty} \sum_{k=0}^\infty 1_{A\cap(V_n=(k+1)/2^n)}
1_{[(k+1)/2^n,\infty)}(t),
\end{align*}
which is $\sO'$ measruable.

(2) As long as $a+(1/n)<b$, the processes $1_C(\omega) 1_{(b-(1/n), c-(1/n)]}(t)$
are left continuous, bounded,  and adapted, hence predictable. 
The process $1_C(\omega)1_{[b,c)}(t)$
is the limit of these processes as $n\to \infty$, so is predictable.
On the other hand, if $X_t$ is a bounded adapted left continuous process,
it can be approximated by 
$$\sum_{k=1}^\infty X_{(k-1)/2^n}(\omega) 1_{(k/2^n, (k+1)/2^n]}(t).$$
Each summand can be approximated by linear combinations of processes
of the form
$1_C(\omega) 1_{(b,c]}(t),$
where $C\in \sF_a$ and $a<b<c$. 
Finally, $1_C 1_{(b,c]}$ is the limit of
$1_C(\omega)1_{[b+(1/n), c+(1/n))}(t)$ as $n\to \infty$.
\qed

A consequence of this lemma is that $\sP\subset \sO$. Since $\sO$ is generated
by the class of right continuous processes and right continuous processes
are progressively measurable, we have from Theorem \ref{GT1} that
the debut of a predictable or  optional set is a stopping time. 

Fix $t$ and define
$$\sO(t)=\{A\cap ([0,t]\times \Omega): A\in \sO\}.$$
Let $\ol \sK^0(t)$ be the collection of subsets of $\sO(t)$               
of the form $K\times C$, where $K$ is a compact subset of $[0,t]$ and
$C\in \sF_a$ with $a\le \inf \{s:s\in K\}$. Let $\ol \sK(t)$ be the collection of finite
unions of sets in $\ol \sK^0(t)$ and $\ol \sK_\delta(t)$ the collection of countable
intersections of sets in $\ol \sK(t)$. 
Define $\ol \sI^0(t)$ to be the collection of sets of the form $[a,b)\times C$,
where $a<b\le t$ and $C\in \sF_a$, 
and let $\ol \sI(t)$ be  the collection of finite unions of sets in $\ol \sI^0(t)$.

The proof of the following proposition  is almost identical to the 
proof of Theorem \ref{GT1}. Because the debut of optional sets is now
known to be a stopping time, it is not nececessary to work with $\P^*$.

\bep\label{SSP1} Suppose $A\in \sO(t)$. Then given $\eps>0$, there
exists $B\in \ol \sK_\delta(t)$ such that $\P(\pi(A)\setminus \pi(B))<\eps$.
\eep

We now prove Theorem \ref{SST1}.

\longproof{ of Theorem \ref{SST1}}
If $E$ is an optional set, choose $t$ large enough so that
if $A_t=E\cap ([0,t]\times \Omega)$, then $\P(\pi(A_t))>\P(\pi(E))-\eps/2$. 
This is possible because $A_t\uparrow E$ and so $\pi(A_t)\uparrow \pi(E)$.
With this value of $t$, choose
$B\in \ol \sK_\delta(t)$ such that $B\subset A_t$ and $\P(\pi(B))>\P(\pi(A_t))-\eps/2$. 
We will show $[D_B,D_B]\subset B$. Since $(D_B<\infty)=\pi([D_B,D_B])=\pi(B)$, we have
$[D_B,D_B]\subset E$
and $$\P(\pi(E))< \P(\pi(A_t))+\eps/2
<\P(\pi(B))+\eps=\P(\pi([D_B,D_B]))+\eps.$$

By the argument of the proof of Lemma \ref{GL1}, $S(B)(\omega)$ is a 
compact set if $B\in \ol \sK_\delta(t)$. Therefore $D_B(\omega)
=\inf\{s:s\in S(B)(\omega)\}$ is in $S(B)(\omega)$, which implies
$[D_B,D_B]\subset B$.
\qed

To prove Theorem \ref{SST2} we follow along the same lines. Define
$$\sP(t)=\{A\cap ([0,t]\times \Omega): A\in \sP\}$$
and define
$\wt \sK^0(t)$ to be the collection of subsets of $\sP(t)$  of the form 
$K\times C$, where $K$ is a compact subset of $[0,t]$ and $C\in \sF_a$
with $a<\inf\{s:s\in K\}$, let $\wt \sK(t)$ be the collection of finite
unions of sets in $\wt \sK^0(t)$, and $\wt \sK_\delta(t)$ the collection of
countable intersections of sets in $\wt \sK(t)$. Define $\wt \sI^0(t)$ to
be the collection of sets of the form $[b,c)\times C$, where $C\in \sF_a$ and
$a<b<c\le t$, and let $\wt \sI(t)$ be the collection of finite unions of sets in
$\wt \sI^0(t)$. 
Following  the proof of Theorem \ref{SST1}, we will be done once we show 
$D_B$ is a predictable stopping time when $B\in \wt \sK_\delta(t)$.

\longproof{ of Theorem \ref{SST2}} Fix $t$. 
Suppose $B\in \wt \sK^0(t)$ is of the form $B=K\times C$ with $C\in \sF_a$ and $a<b=
\inf \{s:s\in K\}$.
Note that this implies $b>0$.
Then  $D_B$ equals $b$ if $\omega\in C$ and equals infinity otherwise. 
As long as $a+(1/m)<b$, we see that $D_A$ is predicted by the stopping times
$V_m$, where $V_m$ equals $b-(1/m)$ if $\omega\in C$ and equals $m$
otherwise.  Note also that $[D_B,D_B]\subset B$. If $B=\cup_{i=1}^m B_i$ with $B_i\in \wt \sK^0(t)$, then
$D_B=D_{B_1}\land \cdots \land D_{B_m}$, and it is easy to see that $D_B$
is predictable because each $D_{B_i}$ is, and also that $[D_B,D_B]\subset B$.

Now let $B\in\wt \sK_\delta(t)$ with $B_n\downarrow B$ and $B_n\in \wt \sK(t)$. We have $D_{B_n}\uparrow$, and the
limit, which we call $T$, will be a stopping time.
Since $B\subset B_n$, then $D_{B_n}\leq D_B$, and therefore $T\leq D_B$.
Each $D_{B_n}$ is a predictable stopping time.  
Let $R_{nm}$ be stopping times predicting $D_{B_n}$ and choose $m_n$ large so
that 
$$\P(R_{nm_n}+2^{-n}<D_{B_n}<\infty)< 2^{-n},
\qq \P(R_{nm_n}<n, D_{B_n}=\infty)<2^{-n}.$$ 
By the Borel-Cantelli lemma,
$$\P(\sup_n R_{nm_n}<T<\infty)=0 \qq \mbox{and} \qq
\P(\sup_n R_{nm_n}<T=\infty)=0,$$ so if we set $Q_n=n\land(R_{1m_1}\lor \cdots \lor
R_{nm_n})$, we see that $\{Q_n\}$ is a sequence
of stopping times predicting $T$, except for a set of probability
zero. 
Hence $T$ is a predictable stopping time.

If $n>m$, then $[D_{B_n},D_{B_n}]\subset B_n\subset B_m$. Since $S(B_m)(\omega)$
is a closed subset of $t$, the facts that $D_{B_n}(\omega)\in S(B_m)(\omega)$ for $n>m$ and
$D_{B_n}(\omega)\to T(\omega)$ for each $\omega$ shows that $T(\omega)\in S(B_m)(\omega)$ for each $\omega$. Thus
$[T,T]\subset B_m$. This is true for all $m$, so $[T,T]\subset B$.
In particular, $T\geq D_B$, so $T=D_B$. Therefore $\pi(B)=(D_B<\infty)
=\pi([T,T])$.

This and the argument of the first paragraph of the proof of Theorem \ref{SST1} proves
Theorem \ref{SST2}.
\qed

\ni{\bf Acknowledgement.} I would like to thank the referee
for valuable suggestions.

%ZZZZZZZZZZZZZZZZZZZZZZZZZZZZZZZZ

\medskip

\ni {\bf Richard F. Bass}\\
Department of Mathematics\\
University of
Connecticut \\
Storrs, CT 06269-3009, USA\\
{\it bass@math.uconn.edu}

\end{document}